\documentclass[reqno]{amsart}

\usepackage{amssymb}
\usepackage{algorithm}
\usepackage{algpseudocode}
\usepackage{graphicx}
\usepackage{subfigure}
\usepackage{epsfig}
\usepackage{url}
\usepackage{cite}

\AtBeginDocument{{\noindent\small
This is a preprint of a paper whose final and definite form will be published
in\\ Journal of Mathematical Analysis, ISSN: 2217-3412, Volume 7, Issue 1 (2016).}
\vspace{9mm}}

\begin{document}

\title[Dengue disease: a multiobjective viewpoint]{Dengue disease:
a multiobjective viewpoint}

\author[R. Denysiuk et al.]{Roman Denysiuk,
Helena Sofia Rodrigues, M. Teresa T. Monteiro,\\
Lino Costa, Isabel Esp\'{i}rito Santo,
Delfim F. M. Torres}

\address{Roman Denysiuk \newline
\indent Algoritmi R\&D Center, University of Minho,
Campus de Gualtar,\newline
\indent 4710-057 Braga, Portugal}
\email{roman.denysiuk@algoritmi.uminho.pt}

\address{Helena Sofia Rodrigues \newline
\indent Algoritmi R\&D Center, University of Minho, Campus de Gualtar,\newline
\indent 4710-057 Braga, Portugal;\newline
\indent School of Business Studies, Viana do Castelo Polytechnic Institute,\newline
\indent Avenida Miguel Dantas, 4930-678 Valen\c{c}a, Portugal;\newline
\indent Center for Research and Development in Mathematics and Applications (CIDMA),\newline
\indent Department of Mathematics, University of Aveiro, 3810-193 Aveiro, Portugal}
\email{sofiarodrigues@esce.ipvc.pt}

\address{M. Teresa T. Monteiro \newline
\indent Department of Production and Systems Engineering, University of Minho,\newline
\indent Campus de Gualtar, 4710-057 Braga, Portugal}
\email{tm@dps.uminho.pt}

\address{Lino Costa \newline
\indent Department of Production and Systems Engineering, University of Minho,\newline
\indent Campus de Gualtar, 4710-057 Braga, Portugal}
\email{lac@dps.uminho.pt}

\address{Isabel Esp\'{i}rito Santo\newline
\indent Department of Production and Systems Engineering, University of Minho,\newline
\indent Campus de Gualtar, 4710-057 Braga, Portugal}
\email{iapinho@dps.uminho.pt}

\address{Delfim F. M. Torres \newline
\indent Center for Research and Development in Mathematics and Applications (CIDMA),\newline
\indent Department of Mathematics, University of Aveiro, 3810-193 Aveiro, Portugal}
\email{delfim@ua.pt}


\thanks{Submitted May 20, 2015. Revised and Accepted Dec 02, 2015.}

\subjclass[2010]{90C29, 92D30, 93A30}

\keywords{Dengue disease; mathematical modelling;
evolutionary multiobjective optimization}


\begin{abstract}
During the last decades, the global prevalence of dengue progressed dramatically.
It is a disease that is now endemic in more than one hundred countries of Africa,
America, Asia, and the Western Pacific. In this paper, we present a mathematical
model for the dengue disease transmission described by a system of ordinary
differential equations and propose a multiobjective approach to find the most
effective ways of controlling the disease. We use evolutionary multiobjective
optimization (EMO) algorithms to solve the resulting optimization problem,
providing the performance comparison of different algorithms. The obtained results
show that the multiobjective approach is an effective tool to solve the problem,
giving higher quality and wider range of solutions compared to the traditional
technique. The obtained trade-offs provide a valuable information about the
dynamics of infection transmissions and can be used as an input in the process
of planning the intervention measures by the health authorities. Additionally,
a suggested hybrid EMO algorithm produces highly superior performance compared
to five other state-of-the-art EMO algorithms, being indispensable
to efficiently optimize the proposed model.
\end{abstract}

\maketitle


\numberwithin{equation}{section}
\newtheorem{theorem}{Theorem}[section]
\newtheorem{lemma}[theorem]{Lemma}
\newtheorem{proposition}[theorem]{Proposition}
\newtheorem{corollary}[theorem]{Corollary}
\newtheorem*{remark}{Remark}


\section{Introduction}

Epidemiology -- the study of patterns of diseases including those which are
noncommunicable infections in population -- has become more relevant and
indispensable in the development of new models and explanations for the outbreaks,
namely due to their propagation and causes. In epidemiology, an infection
is said to be endemic in a population when it is maintained in the population
without the need for external inputs. An epidemic occurs when new cases
of a certain disease appears, in a given human population during a given period,
and then essentially disappears.

Mathematical modeling is critical to understand how epidemiological diseases
spread. It can help to explain the nature and dynamics of infection transmissions
and can be used to devise effective strategies for fighting them. When formulating
a model for a particular disease, we should make a trade-off between simple models
-- that omit several details and generally are used for specific situations
in a short time, but have the disadvantage of possibly being naive and unrealistic
-- and more complex models, with more details and more realistic, but generally
more difficult to solve or could contain parameters which their estimates
cannot be obtained.

Dengue is a vector-borne disease transmitted from an infected human to a female
\emph{Aedes} mosquito by a bite. Then, the mosquito, that needs regular meals
of blood to feed their eggs, bites a potentially healthy human and transmits
the disease, turning it into a cycle.

There are four distinct, but closely related, viruses that cause dengue.
The four serotypes, named DEN-1 to DEN-4, belong to the \emph{Flavivirus} family,
but they are antigenically distinct. Recovery from infection by one virus provides
lifelong immunity against that virus but provides only partial and transient
protection against subsequent infection by the other three viruses. There are
strong evidences that a sequential infection increases the risk of developing
dengue hemorrhagic fever.

The spread of dengue is attributed to the geographic expansion of the mosquitoes
responsible for the disease: \emph{Aedes aegypti} and \emph{Aedes albopictus}.
The \emph{Aedes aegypti} mosquito is a tropical and subtropical species widely
distributed around the world, mostly between latitudes $35^{o}$N and 35$^o$S.
In urban areas, \emph{Aedes} mosquitoes breed on water collections in artificial
containers such as cans, plastic cups, used tires, broken bottles and flower pots.
Due to its high interaction with humans and its urban behavior, the
\emph{Aedes aegypti} mosquito is considered the major responsible
for the dengue transmission.

\begin{figure}
\center
\includegraphics[scale=0.55]{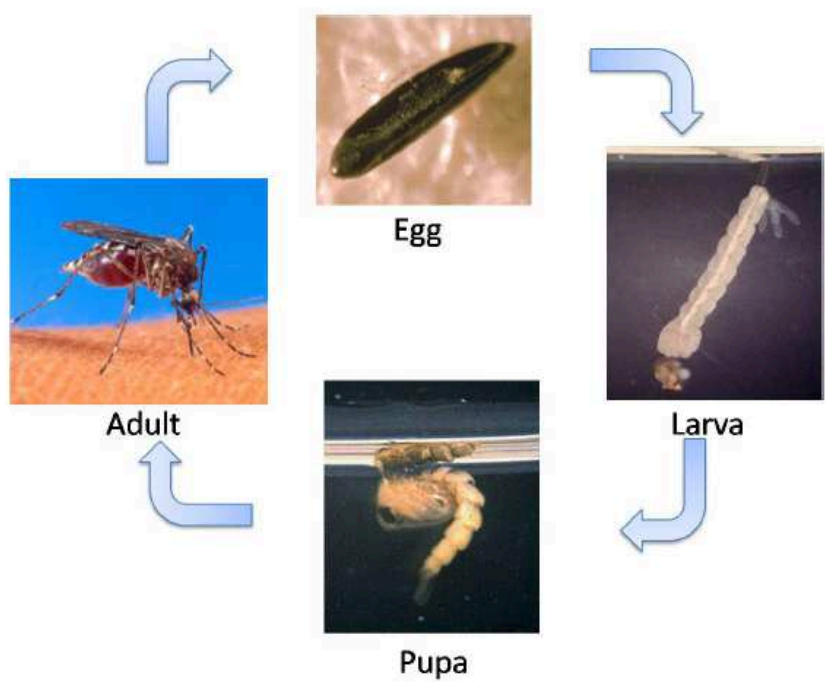}
\caption[Life cycle of \emph{Aedes aegypti}.]{Life cycle of \emph{Aedes aegypti}.}
\label{app:fig:lifecycle}
\end{figure}

The life cycle of a mosquito has four distinct stages: egg, larva, pupa and adult,
as it is possible to see in Figure~\ref{app:fig:lifecycle}. In the case of
\emph{Aedes aegypti}, the first three stages take place in, or near, the water,
whereas the air is the medium for the adult stage~\cite{Otero2008}.

It is very difficult to control or eliminate \emph{Aedes aegypti} mosquitoes
due to their resiliency, fast adaptation to changes in the environment
and their ability to rapidly bounce back to initial numbers after disturbances
resulting from natural phenomena (e.g., droughts) or human interventions
(e.g., control measures).

Primary prevention of dengue resides mainly in mosquito control. There are two
primary methods: larval control and adult mosquito control, depending
on the intended target. \emph{Larvicide} treatment is done through long-lasting
chemical in order to kill larvae and preferably have World Health Organization
clearance for use in drinking water~\cite{Derouich2003}. \emph{Adulticides} is
the most common measure. Its application can have a powerful impact on the
abundance of adult mosquito vector. However, the efficacy is often constrained
by the difficulty in achieving sufficiently high coverage
of resting surfaces~\cite{Devine2009}.

In this paper, our first motivation is to analyse a mathematical model for the
dengue disease transmission, including the control parameter representing
measures to fight the disease. Since there are different goals that can be sought
dealing with the dengue epidemic, we suggest a multiobjective approach to find
the optimal control in the proposed model. Such approach better reflects the
nature of an underlying decision-making problem, hence,
can provide more interesting results.

Without loss of generality, a multiobjective optimization (MO) problem with $m$
objectives and $n$ decision variables can be mathematically formulated as follows:
\begin{equation}
\begin{array}{rl}
\text{minimize:} & \boldsymbol{f}(\boldsymbol{x})
=(f_1(\boldsymbol{x}),f_2(\boldsymbol{x}),\ldots,f_m(\boldsymbol{x}))^{\text{T}}\\
\text{subject to:} & \boldsymbol{x} \in \Omega,
\end{array}
\end{equation}
where $\boldsymbol{x}$ is the decision vector, $\Omega \subseteq \mathbb{R}^n$
is the feasible decision space, and $\boldsymbol{f}(\boldsymbol{x})$ is the
objective vector defined in the objective space $\mathbb{R}^m$.

When several objectives are optimized at the same time, the search space becomes
partially ordered. In such scenario, solutions are compared on the basis of the
Pareto dominance. For two solutions $\boldsymbol{a}$ and $\boldsymbol{b}$ from
$\Omega$, a solution $\boldsymbol{a}$ is said to dominate a solution
$\boldsymbol{b}$ (denoted by $\boldsymbol{a} \prec \boldsymbol{b}$) if:
\begin{equation*}
\forall i \in \{ 1, \ldots ,m\} :f_i (\boldsymbol{a}) \le f_i (\boldsymbol{b})
\, \wedge \, \exists j \in \{ 1, \ldots ,m\}
:f_j (\boldsymbol{a}) < f_j (\boldsymbol{b}).
\end{equation*}
Since solutions are compared against different objectives, there is no longer
a single optimal solution but a set of optimal solutions, generally known as
the Pareto optimal set. This set contains equally important solutions representing
different trade-offs between the given objectives and can be defined as:
\begin{equation*}
\mathcal{PS} = \{ \boldsymbol{x} \in \Omega \, | \, \nexists \boldsymbol{y}
\in \Omega : \boldsymbol{y} \prec \boldsymbol{x} \}.
\end{equation*}
Approximating the Pareto optimal set is the main goal in multiobjective optimization.

Evolutionary multiobjective optimization (EMO) algorithms have become popular
and successful in solving MO problems. Like their single-objective counterparts,
EMO algorithms mimic the principles of the natural selection and evolution.
EMO algorithms are population-based optimization techniques. This feature makes
them especially attractive to dealing with MO problems, allowing to approximate
the Pareto optimal set in a single simulation run. Since the first implementations
appeared in the mid 1980s, there has been a growing interest to developing
efficient EMO algorithms. This resulted in a plenty of proposals~\cite{CoelloBook,DebBook}.
All these algorithms proved their viability in a number of comparative studies
and were successfully used in many real-world applications~\cite{ZhouQu2011}.
However, it is well-known that overall successful and efficient general solvers
do not exist. Statements about the optimization algorithms' performance must
be qualified with regard to the ``no free lunch'' theorem for
optimization~\cite{WolpertMacready1997}. Its simplest interpretation is that
a general-purpose universal optimization strategy is theoretically impossible,
and the only way one strategy can outperform another is if it is specialized
to the specific problem under consideration~\cite{HoPepyne2002}. This motivates
to study strengths and weaknesses of existing techniques, contributing to the
development of new efficient optimization methods. Therefore, additionally to
solving the problem resulting from the model for the dengue disease, we perform
a comparative study of different EMO algorithms on the proposed model. Motivated
by the promising results presented in~\cite{Denysiuk6}, we suggest an improved
descent direction-guided multiobjective algorithm (DDMOA2) as the base algorithm
for solving the herein formulated MO problem. DDMOA2 is a hybrid multiobjective
algorithm inspired by its predecessor, which proved to be highly competitive
and often superior compared to the existing approaches~\cite{Denysiuk5}. For
the comparison, we select the most representative state-of-the-art EMO algorithms
with popular evolutionary operators to produce offspring and different algorithmic
frameworks. Thus, our second motivation in this paper is to promote a hybrid
methodology in the design of EMO algorithms. Such methods do not only increase
the variety of existing approaches, they also represent a significant contribution
to the optimization, often becoming indispensable to successfully
solve real-world problems.

The remainder of this paper is organized as follows. Section~\ref{sec:2}
describes a mathematical model for the dengue disease transmission and
formulates a multiobjective optimization problem to find the optimal control
in the presented model. Section~\ref{sec:3} gives a detailed description of
DDMOA2 and provides an overview of five other state-of-the-art EMO algorithms.
Section~\ref{sec:4} offers the performance comparison of all the algorithms,
analysing the difficulties faced by the solvers. It also discusses different
perspectives of the dengue disease extracted from the obtained trade-off
solutions, and compares the proposed approach with the traditional one.
Section~\ref{sec:5} concludes the work.


\section{Mathematical Model}
\label{sec:2}

This section introduces a mathematical model for the dengue disease transmission
based on a system of ordinary differential equations and the real data of a
dengue disease outbreak that occurred in the Cape Verde archipelago in
2009 \cite{ministerio,RMT2013a}. After describing the model, a multiobjective
approach is proposed to find the most effective ways of controlling the disease.


\subsection{ODE SEIR+ASEI Model with Insecticide Control}

The model consists of eight mutually-exclusive compartments representing
the human and vector dynamics. It also includes a control parameter,
an adulticide spray, as a measure to fight the disease.

The notation used in the mathematical model includes
four epidemiological states for humans:
\begin{quote}
\begin{tabular}{lcl}
$S_h(t)$ & -- & susceptible; \\
$E_h(t)$ & -- & exposed; \\
$I_h(t)$ & -- & infected; \\
$R_h(t)$ & -- & resistant.
\end{tabular}
\end{quote}
It is assumed that the total human population $(N_h)$
is constant, so, $N_h=S_h+E_h+I_h+R_h$.

There are also four other state variables related to the mosquitoes:
\begin{quote}
\begin{tabular}{lcl}
$A_m(t)$ & -- & aquatic phase; \\
$S_m(t)$ & -- & susceptible; \\
$E_m(t)$ & -- & exposed; \\
$I_m(t)$ & -- & infected. \\
\end{tabular}
\end{quote}
Similarly, it is assumed that the total adult mosquito population
is constant, which means $N_m=S_m+E_m+I_m$.

The model includes a control variable, which represents the amount
of insecticide that is continuously applied during a considered period,
as a measure to fight the disease:
\begin{quote}
\begin{tabular}{lcl}
$c(t)$ & -- & level of insecticide campaigns.
\end{tabular}
\end{quote}
The control variable is an adimensional value that is considered
in relative terms varying from 0 to 1.

In the following, for the sake of simplicity, the independent variable $t$ is
omitted when writing the dependent variables (for instance,
$S_h$ is used instead of $S_h (t)$).

The parameters necessary to completely describe the model are presented in
Table~\ref{tab:param}. This set of parameters includes the real data related
to the outbreak of dengue disease occurred in Cape Verde in 2009.

\begin{table}[ht!]
\centering
\caption{Model parameters.}\label{tab:param}
\begin{tabular}{|c|c|c|}
\hline
Parameter & Meaning& Value \\ \hline
$N_h$ & total population & 480000 \\[0.2cm]
$B$ & average daily bites (per mosquito per day)& 1 \\[0.2cm]
$\beta_{mh}$ & transmission probability from $I_m$ (per bite)& 0.375  \\[0.2cm]
$\beta_{hm}$ & transmission probability from $I_h$ (per bite)& 0.375  \\[0.2cm]
$\mu_{h}$ & average human lifespan (in days)& $1/(71\times365)$  \\[0.2cm]
$\eta_{h}$ & mean viremic period (in days)& 1/3 \\[0.2cm]
$\mu_{m}$ & average lifespan of adult mosquitoes (in days)& 1/11  \\[0.2cm]
$\varphi$ & number of eggs at each deposit per capita (per day)& 6  \\[0.2cm]
$\mu_{A}$ & natural mortality of larvae (per day)& 1/4  \\[0.2cm]
$\eta_A$ & rate of maturing from larvae to adult (per day)& 0.08  \\[0.2cm]
$\eta_m$ & extrinsic incubation period (in days)& 1/11   \\[0.2cm]
$\nu_h$ & intrinsic incubation period (in days)& 1/4  \\[0.2cm]
$m$ & number of female mosquitoes per human & 6 \\[0.2cm]
$k$ & number of larvae per human & 3\\ \hline
\end{tabular}
\end{table}

Furthermore, in order to obtain a numerically stable problem,
all the state variables are normalized as follows:
\begin{center}
\begin{tabular}{llll}
$\displaystyle s_h=\frac{S_h}{N_h},$ & $\displaystyle e_h
=\frac{E_h}{N_h},$ & $\displaystyle i_h=\frac{I_h}{N_h},$ & $\displaystyle r_h=\frac{R_h}{N_h},$\\
& & & \\
$\displaystyle a_m=\frac{A_h}{kN_h},$ & $\displaystyle s_m=\frac{S_m}{mN_h},$
& $\displaystyle e_m=\frac{E_m}{mN_h},$ & $\displaystyle i_m=\frac{I_m}{mN_h}.$\\
\end{tabular}
\end{center}
Thus, the dengue epidemic is modeled by the following nonlinear
time-varying state equations~\cite{RodriguesMonteiro2012}:
\begin{equation}
\label{app:dengue:eq:1}
\begin{tabular}{l}
$\left\{
\begin{array}{l}
\displaystyle\frac{ds_h}{dt} = \mu_h - (B\beta_{mh}m i_m +\mu_h) s_h\\
\displaystyle\frac{de_h}{dt} = B\beta_{mh}m i_m s_h - (\nu_h + \mu_h )e_h\\
\displaystyle\frac{di_h}{dt} = \nu_h e_h -(\eta_h  +\mu_h) i_h\\
\displaystyle\frac{dr_h}{dt} = \eta_h i_h - \mu_h r_h\\
\displaystyle\frac{da_m}{dt} = \varphi \frac{m}{k}(1-a_m)(s_m+e_m+i_m)-(\eta_A+\mu_A) a_m\\
\displaystyle\frac{ds_m}{dt} = \eta_A \frac{k}{m}a_m-(B \beta_{hm}i_h+\mu_m) s_m-c s_m\\
\displaystyle\frac{de_m}{dt} = B \beta_{hm}i_h s_m-(\mu_m + \eta_m) e_m-c e_m\\
\displaystyle\frac{di_m}{dt} = \eta_m e_m -\mu_m i_m - c i_m\\
\end{array}
\right. $
\end{tabular}
\end{equation}
\noindent with the initial conditions
\begin{center}
\begin{tabular}{llll}
$s_h(0)=0.99865,$ & $e_h(0)=0.00035,$ & $i_h(0)=0.001,$ &
$r_h(0)=0,$ \\
$a_m(0)=1,$ & $s_{m}(0)=1,$ &
$e_m(0)=0,$ & $i_m(0)=0.$
\end{tabular}
\end{center}
Since any mathematical model is an abstraction of a complex natural system,
additional assumptions are made to make the model mathematically treatable.
This is also the case for the above epidemiological model,
which comprises the following assumptions:
\begin{itemize}
\item the total human population ($N_h$) is constant;
\item there is no immigration of infected individuals into the human population;
\item the population is homogeneous, which means that every individual
of a compartment is homogeneously mixed with other individuals;
\item the coefficient of transmission of the disease is fixed
and does not vary seasonally;
\item both human and mosquitoes are assumed to be born susceptible,
i.e., there is no natural protection;
\item there is no resistant phase for mosquitoes, due to their short lifetime.
\end{itemize}

The mathematical model for the dengue disease given by~\eqref{app:dengue:eq:1}
includes the control parameter, which represents a measure that can be taken
by authorities in response to the disease. Analysing the model, a natural
question arises -- What is the most effective way of applying the control?
A problem of finding a control law for a given system is commonly formulated
and solved using Optimal Control Theory, where a control problem includes
a cost functional that is a function of state and control variables. Thus,
the objective functional $J$, considering the costs of infected humans
and costs with insecticide, can formulated as follows~\cite{RodriguesMonteiro2013}:
\begin{equation}
\label{functional}
\text{minimize } J=\int_{0}^{T}\left[ \gamma_D i_h(t)^2 + \gamma_S c(t)^2\right]dt,
\end{equation}
where $\gamma_D$ and $\gamma_S$ are positive constants representing the costs
weights of infected individuals and spraying campaigns, respectively.
The optimal control can be derived using Pontryagin's
Maximum Principle~\cite{pontryagin62}.


\subsection{Multiobjective Approach}

An approach based on Optimal Control Theory allows to obtain a single optimal
solution. The most straightforward disadvantage of such approach is that only
a limited amount of information about a choice of the optimal strategy can be
presented to the decision maker. On the other hand, when formulating a
multiobjective optimization problem one seeks to simultaneously optimize
several objectives. Solving a multiobjective optimization problem yields
a set of optimal solutions, given two or more conflicting objectives. This set
provides a whole range of optimal strategies, being crucial to an effective
decision-making process. Thus, we propose a multiobjective approach to find
optimal strategies for applying insecticide. Our approach is based
on simultaneously optimizing the cost due to infected population and the cost
associated with the insecticide. The problem is stated as follows:
\begin{equation}
\label{app:dengue:eq:2}
\begin{array}{rl}
\text{minimize:}   & f_1(c(t)) = \int_{0}^{T} i_h(t) \, dt \\
& f_2(c(t)) = \int_{0}^{T} c(t) \, dt \\
\text{subject to:} & \eqref{app:dengue:eq:1}
\end{array}
\end{equation}
where $T$ is the period of time, $f_1$ and $f_2$ represent the total cost
incurred in the form of infected population and the total cost of applying
insecticide for the period $T$, respectively.


\section{Algorithms and Experimental Design}
\label{sec:3}

This section describes EMO algorithms and the experimental design used to find
the optimal control in the mathematical model for the dengue disease transmission
given by~\eqref{app:dengue:eq:2}. First, we give a detailed description of the
second version of descent directions-based multiobjective algorithm. It follows
a brief discussion of five other state-of-the-art EMO algorithms. Finally,
the general experimental setting is provided for all the algorithms.


\subsection{DDMOA2}

\begin{algorithm}[b!]
\caption{DDMOA2}\label{ddmoa2:alg:1}
\small
\begin{algorithmic}[1]
\State $g\gets 0$;
\State $\texttt{initialization:}$ $P^{(g)}$,
$W=\{\boldsymbol{w}_1,\ldots,\boldsymbol{w}_{\mu} \}$;
\Repeat
\State $P^{(g)} \gets \texttt{leaderSelection}(P^{(g)})$;
\State $P^{(g)} \gets \texttt{updateSearchMatrix}(P^{(g)})$;
\State $P^{(g)} \gets  \texttt{updateStepSize}(P^{(g)})$;
\State $P^{(g)} \gets \texttt{parentSelection}(P^{(g)})$;
\State $P^{(g)} \gets \texttt{mutation}(P^{(g)})$;
\State $P^{(g+1)} \gets \texttt{environmentalSelection}(P^{(g)})$;
\State $g \gets g + 1$;
\Until{the stopping criterion is met}
\State  output: $P^{(g)}$;
\end{algorithmic}
\end{algorithm}

Descent directions-based multiobjective algorithm (DDMOA2) described in this
paper is a recent hybrid multiobjective optimization algorithm proposed
by Denysiuk et al.~\cite{Denysiuk6}. The main loop of DDMOA2 is given by
Algorithm~\ref{ddmoa2:alg:1}. In each generation, the selection of leaders
and the adaptation of the strategy parameters of all population members
are followed by the successive application of parent selection, mutation,
and environmental selection.

In DDMOA2, an individual $a_i$ ($i \in \{1,\ldots,\mu\}$) in the current
population $P^{(g)}$ in generation $g$ is a tuple of the form
$[\boldsymbol{x}_i,\delta_i,\boldsymbol{S}_i,\sigma_i]$, where
$\boldsymbol{x}_i \in \mathbb{R}^n$ is the decision vector, $\delta_i>0$
is the step size used for local search, $\boldsymbol{S}_i
\in \mathbb{R}^{n\times 2}$ is the search matrix, and $\sigma_i>0$
is the step size used for reproduction.

In the following, the components of DDMOA2 are discussed in more detail.


\subsubsection{Initialization Procedure}
\label{ddmoa2:subsec:1}

The algorithm starts by generating a set of weight vectors
$W=\{\boldsymbol{w}_1,\ldots,\boldsymbol{w}_{\mu} \}$ and
initializing the population of size $\mu$ using Latin hypercube
sampling~\cite{Loh1996}. The strategy parameters of each population
member are initialized taking default values. The search matrix
$\boldsymbol{S}^{(0)}$ is initialized by simply generating
a zero matrix of size $n\times 2$.


\subsubsection{Leader Selection Procedure}

Each generation of DDMOA2 is started by selecting leaders of the current
population. A leader is a population member that performs the best on at
least one weight vector. Thus, leaders are selected as follows. First,
the objective values of all individuals in the population are normalized:
\begin{equation}
\label{ddmoa2:eq:2}
\overline{f}_i=\frac{f_i-f_i^{\text{min}}}{f_i^{\text{max}}-f_i^{\text{min}}},
\quad \forall i \in \{1,\ldots,m\},
\end{equation}
where $f_i^{\text{min}}$ and $f_i^{\text{max}}$ are the minimum and maximum
values of the $i$th objective in the current population, respectively, and
$\overline{f}_i \in [0,1], \forall i \in \{1,\ldots,m\}$ is the normalized
objective value. For each weight vector, the fitness of each population member
is calculated on a given weight vector using the weighted Chebyshev method,
which after normalization of objectives can be defined as:
\begin{equation}
\label{ddmoa2:eq:3}
f_{\text{fitness}}=\max \limits_{1\leq i\leq m}
\,\, \{\,w_i \, \overline{f}_i(\boldsymbol{x})\,\},
\end{equation}
where $f_{\text{fitness}}$ is the fitness of the population member
$\boldsymbol{x}$ on the weight vector $\boldsymbol{w}$. An individual having
the best fitness on a given weight vector is a leader. It should be noted that
one leader can have the best fitness value on several weight vectors.


\subsubsection{Update Search Matrix Procedure}

After selecting leader individuals, the search matrices of all individuals
in the current population are updated in the \texttt{update\-Search\-Matrix}
procedure. In DDMOA2, descent directions are calculated only for two randomly
chosen objectives, regardless of the dimensionality of the objective space.
Thus, the search matrix of each population member contains two columns that
store descent directions for these randomly chosen objectives.

In the beginning of the \texttt{updateSearchMatrix} procedure, two objectives
are chosen at random, and only leaders of the current population are considered
while all the other individuals are temporarily discarded. Thereafter, for the
first chosen objective, the resulting population is sorted in ascending order
and partitioned into $\alpha$ equal parts. Thus, $\alpha$ subpopulations are
defined in order to promote different reference points for the computation
of descent directions. It follows that in each subpopulation, a representative
individual $a_{\text{r}}$ is selected. A representative of the subpopulation
is a solution with the smallest value of the corresponding objective function
among other solutions in the subpopulation and $\delta>\delta_{\text{tol}}$.
Thus, if the solution with the smallest value of the corresponding objective
has $\delta \leq \delta_{\text{tol}}$ then the solution with the second smallest
value is selected as representative and so on. After that, a descent direction
for the corresponding objective function is computed for the representative
using coordinate search \cite{Torczon1997}. It should be noted that a descent
direction $\boldsymbol{s}$ for subpopulation representative
$\boldsymbol{x}_{\text{r}}$ is accepted when trial solution
$\boldsymbol{x}_{\text{r}}+\boldsymbol{s}$ is nondominated with respect
to the current population. During coordinate search, the step size $\delta$
of the representative is reduced if no decrease in the objective function
value is found. Each time a trial solution is calculated, this solution
is compared with the current population. If this trial solution has a smaller
value for at least one objective compared with each member of the current
population (i.e. the trial is nondominated with respect to the current population)
then it is added to the population, assuming the default values of the strategy
parameters. When a descent direction $\boldsymbol{s}_{\text{r}}$ for the
subpopulation representative is found, descent directions for all other
subpopulation members are computed as follows:
\begin{equation}
\label{ddmoa2:eq:4}
\boldsymbol{s}_i=\boldsymbol{x}_{\text{r}}
-\boldsymbol{x}_i+\boldsymbol{s}_{\text{r}},
\end{equation}
where $\boldsymbol{s}_i$ is the descent direction for the $i$th
subpopulation member, $\boldsymbol{x}_i$ is the decision vector of the
$i$th subpopulation member, $\boldsymbol{x}_{\text{r}}$ is the subpopulation
representative, and $\boldsymbol{s}_{\text{r}}$ is the descent direction for
the representative. The calculated descent directions are stored in the first
column of the search matrices of the corresponding individuals. Thereafter,
the same procedure for finding descent directions is performed for the second
chosen objective and the results are stored in the second column
of the search matrices.

After the search matrices of leaders of the current population are updated,
the search matrices of all other population members need to be updated. For
this purpose, a simple stochastic procedure is used. For each non-leader
individual, a leader is randomly chosen and this leader shares its search
matrix, i.e., non-leader individual's search matrix is equal to the selected
leader's search matrix. At the end of the \texttt{updateSearchMatrix} procedure,
the search matrices of all population members are updated. Since some promising
solutions satisfying the aforementioned conditions are added to the population
during coordinate search, at the end of this procedure, the population size
is usually greater than $\mu$.


\subsubsection{Update Step Size Procedure}

Before generating offspring, the step size of each population member needs
to be updated. However, there is no common rule to update the step size
$\sigma$, but it must be done carefully to ensure convergence to the Pareto set,
starting from a larger value at the beginning and gradually reducing it during
the generations. DDMOA2 uses the following rule for updating the step size
of each population member:
\begin{equation}
\sigma=\max \{ \exp(\tau N(0,1)) \, \sigma_0^{(1-\frac{3\,funEval}{maxEval})},
\delta_{\text{tol}}\},
\end{equation}
where $\tau$ is the learning parameter ($\tau=1/\sqrt{2n}$), $N(0,1)$ is a random
number sampled from the normal distribution with mean $0$ and standard deviation
$1$, $\sigma_0$ is the initial value of the step size, $funEval$ is the current
number of function evaluations, and $maxEval$ is the maximum number of function
evaluations. The idea here is to exponentially reduce the step size depending
on the current number of function evaluations and multiply it by a scaling
factor (borrowed from evolution strategies~\cite{BeyerSchwefel2002}), thereby
obtaining different step size values among the population members.


\subsubsection{Parent Selection Procedure}

In each generation, $\lambda$ offspring individuals are generated by mutating
the correspondent parent. The decision on how many offspring are produced from
each population member is made in the \texttt{parentSelection} procedure. At
the beginning, it is assumed that none offspring is generated by each individual.
Then, binary tournament selection based on scalarizing fitness is performed
to identify which population members will be mutated in order to produce
offspring. This selection process is based on the idea proposed in
\cite{Ishibuchi2006} and fosters promising individuals to produce more offspring.

The procedure starts by normalizing the objective values of all individuals
in the population as defined in~\eqref{ddmoa2:eq:2}. Then, the selection process
is performed in two stages. In the first stage, only leaders of the current
population are considered. For each weight vector, two individuals are randomly
selected and the fitness of the corresponding individual is calculated as defined
in \eqref{ddmoa2:eq:3}. The individual having the smallest fitness value
is a winner and the number of times it is mutated is augmented by one. Thereafter,
all leaders are removed from the population and, for each weight vector, binary
tournament selection is performed on the resulting population in the same way
as it is done for leaders. It should be noted that some population members might
be mutated several times, as a result of this procedure, while others
will not produce any offspring at all.


\subsubsection{Mutation Procedure}

After identifying which population members have to be mutated, offspring are
generated in the \texttt{mutation} procedure. For the corresponding individual,
the mutation is performed as follows:
\begin{equation}
\label{ddmoa2:eq:5}
\boldsymbol{x}' = \boldsymbol{x}+ \sigma \boldsymbol{S} \, \boldsymbol{\nu},
\end{equation}
where $\sigma$ is the step size, $\boldsymbol{S}$ is the search matrix and
$\boldsymbol{\nu}$ is a column vector of random numbers sampled from the uniform
distribution ($\forall i \in \{1,2\}: \nu_i \sim \mathbb{U} (0,1)$). To guarantee
that each new solution $\boldsymbol{x}'=(x_1',\ldots,x_n')^{\text{T}}$ belongs
to $\Omega$, projection is applied to each component of the decision vector:
$\boldsymbol{x}' = \min\{\max\{\boldsymbol{x}',\boldsymbol{l}\},\boldsymbol{u}\}$.
After offspring $\boldsymbol{x}'$ is repaired, it is evaluated and added
to the population.


\subsubsection{Environmental Selection Procedure}

At the end of each generation, $\mu$ fittest individuals are selected from the
enlarged population in the \texttt{environmental\-Se\-lec\-tion} procedure.
DDMOA2 uses the selection mechanism proposed in \cite{Hughes2007}.

First, the objective values of all individuals are normalized as defined
in Equation~\eqref{ddmoa2:eq:2}. Next, the following steps are performed:
\begin{enumerate}
\item matrix $\boldsymbol{M}$ is calculated, which stores metrics for each
population member on each weight vector (for each population member, a metric
on a corresponding weight vector is computed as defined in~\eqref{ddmoa2:eq:3});

\item for each column (weight vector), the minimum and second smallest
metric value are found;

\item for each column, metric values are scaled by the minimum value found,
except for the row which gave the minimum value. This result is scaled
by the second lowest value;

\item for each row (population member), the minimum scaled value is found,
this value represents individual's fitness;

\item the resulting column vector is sorted, and $\mu$ individuals
with the smallest fitness values are selected.
\end{enumerate}

Normalizing objective values allows to cope with differently scaled objectives,
while the weighted Chebyshev method can find optimal solutions in convex
and nonconvex regions of the Pareto front. When the stopping criterion is met,
the algorithm returns its final population.


\subsection{Other EMO Algorithms}

Nondominated sorting genetic algorithm (NSGA-II) was suggested by
Deb et al.~\cite{DebPratap2002}. In NSGA--II, each solution
in the current population is evaluated using the Pareto ranking and
a crowding distance measure. Fitness assignment procedure ranks the
population according to different non-domination levels using the fast
nondominated sorting. In each non-domination level, the crowding distance
is calculated and assigned to each solution. During the mating selection
procedure, the binary tournament selection based on two sorting criteria,
namely, the non-domination rank and crowding distance is performed to select
a mating pool. When two solutions are selected, the one with the lower
non-domination rank is preferred. If both solutions belong to the same rank,
then the solution with the higher crowding distance is selected. The offspring
population is created by selecting two parents from the mating pool and applying
genetic operators. The real representation of individuals is used, and the
genetic operators are: the SBX crossover and the polynomial mutation. After
the offspring population is generated, it is combined with the parent population
and the fitness assignment procedure is performed. Then, the environmental
selection procedure is used to select the population of the next generation.
Thus, each non-domination level is selected one at a time, starting from the
first level, until no further level can be included without increasing
the population size. In general, the last accepted level cannot be completely
accommodated. In such a case, only solutions having higher values of the crowding
distance are selected to complete the population.

Indicator-based evolutionary algorithm (IBEA) was proposed by Zitzler
and K\"{u}nzli~\cite{ibea}. IBEA calculates the individuals' fitness
by comparing pairs of population members using binary performance indicators.
In this paper, IBEA is used in combination with the additive epsilon indicator,
which subsumes the translations in each dimension of the objective space that
are necessary to create a weakly dominated solution. In the mating selection
procedure, the binary tournament selection with replacement based on the scalar
fitness values is performed to fill a mating pool. In the variation procedure,
two parents are selected from the mating pool and genetic operators, namely,
the SBX crossover and the polynomial mutation, are applied to produce offspring.
After that, the parent and offspring populations are combined. During the
environmental selection procedure used to select the population of the next
generation, individuals having the worst fitness values are iteratively removed
from the population, updating the fitness values of the remaining individuals.

The third version of a generalized differential evolution (GDE3) was proposed
by Kukkonen and Lampinen~\cite{KukkonenLampinen2005}. The proposal extends
the $DE\slash rand\slash 1\slash bin$ variant to solve multiobjective
optimization problems. In each generation, GDE3 goes through population members
and creates a corresponding trial vector using the differential evolution (DE)
operator. When the trial is generated, it is compared with the target vector.
The selection of the solution for the next generation is based on the following
rules. In the case of infeasible vectors, the trial vector is selected if it
weakly dominates the target vector in the constraint violation space, otherwise
the target vector is selected. In the case of the feasible and infeasible vectors,
the feasible vector is selected. If both vectors are feasible, then the trial
is selected if it weakly dominates the target vector in the objective space.
If the target vector dominates the trial vector, then the target vector
is selected. If neither vector dominates each other in the objective space,
then both vectors are kept in the population. After a generation,
the population size may increase. In this case, the population is sorted using
a modified nondominated sorting, which takes into consideration constraints.
As the second sorting criterion, the crowding distance measure is used.
The best population members according to these two sorting criteria
are selected to form the population of the next generation.

Multiobjective evolutionary algorithm based on decomposition (MOEA/D)
was suggested by Li and Zhang~\cite{LiZhang2009}. MOEA/D decomposes
a multiobjective optimization problem into a number of scalar optimization
subproblems and optimizes them simultaneously in a single run instead
of solving a multiobjective problem directly. Neighborhood relations among
the single objective subproblems are defined based on the distances among
their weight vectors. For each weight vector, a certain number of closest
weight vectors are associated, which constitute the neighborhood to this
vector. Each subproblem is optimized by using information from its several
neighboring subproblems. Each individual in the current population is
associated with the corresponding subproblem. In each generation,
MOEA/D goes through each subproblem and generates an offspring using
the DE operator, where the individual associated with the given subproblem
is used as the parent. The other two individuals are selected from its
neighborhood with a certain probability, or from the whole population
otherwise. If a child solution is of high quality, it replaces current
solutions performing worse in its neighborhood or in the entire population.
An extra parameter is introduced to control the maximum number
of solutions replaced by a child solution.

Speed-constrained multiobjective particle swarm optimization (SMPSO) was
proposed by Nebro et al.~\cite{NDG09}. In each flight cycle, SMPSO maintains
a swarm that includes the position, velocity, and individual best position
of the particles. A new position of the particle is calculated by adding
an updated speed to the previous position of the particle. The speed of
each particle is calculated taken into account its best position and the
leader particle taken from the archive containing nondominated solutions.
SMPSO uses a special velocity constriction procedure to control the speed
of the particles. After the velocity for each particle is calculated using
the traditional particle swarm optimization rule, the obtained velocity
is multiplied by a constriction factor. Then, the resulting velocity value
is constrained by bounds calculated for each component of the velocity.
Additionally, the polynomial mutation is applied on each particle to accelerate
the convergence of the swarm. After positions of the particles are updated,
the swarm is evaluated. Thereafter, the memory of each particle and the archive
of leaders are updated. To maintain the archive of bounded size, the crowding
distance measure is used to decide which particles must be retained in the
archive. At the end of the optimization process, the leaders' archive contains
the Pareto set approximation.


\subsection{Experimental Setup}

The constraints presented in~\eqref{app:dengue:eq:1} are numerically integrated
using the fourth-order Runge-Kutta method with 1000 equally spaced time
intervals over the period of 84 days (the Cape Verde outbreak period). This way,
the system of differential equations is solved, and its components become
a part of the objective functions given in problem~\eqref{app:dengue:eq:2},
thereby obtaining an unconstrained problem. The decision variable $c(t)$
is also discretized ending up with a total of 1001 decision variables. Thus,
the feasible decision space is $c \in [0,1]^{1001}$. Furthermore, the integrals
used to determine objective function values in~\eqref{app:dengue:eq:2} are
calculated using the trapezoidal rule. The problem is coded in the
MATLAB\textsuperscript{\textregistered} and Java\textsuperscript{\tiny TM}
programming languages.

The MATLAB\textsuperscript{\textregistered} implementation of DDMOA2 is used,
while the other five EMO algorithms are used within the jMetal framework
\cite{jmetal}. For each algorithm, 30 independent runs are performed with
a population size of $100$, running for $10^5$ function evaluations. For each
algorithm within the jMetal framework, the remaining parameters use the default
settings. The parameter settings employed for DDMOA2 are: the initial step
size for local search $\delta^{(0)}=0.4$, the initial step size for reproduction
$\sigma^{(0)}=5$, the number of subpopulations $\alpha=5$, and the tolerance
for step size $\delta_{\text{tol}}=10^{-3}$.


\section{Discussion}
\label{sec:4}

This section reports the empirical results obtained after the performed
experiments. The performance comparison of the algorithms is carried out,
and the problem difficulties are discussed. The obtained trade-off solutions
are presented, and different scenarios of the dengue disease transmission
are considered. Finally, we provide the comparison of the results obtained
for this problem using the proposed multiobjective approach and the
traditional technique presented in \eqref{functional}.


\subsection{Performance Comparison}

\begin{figure}
\centering
\subfigure[NSGA-II]{\epsfig{file=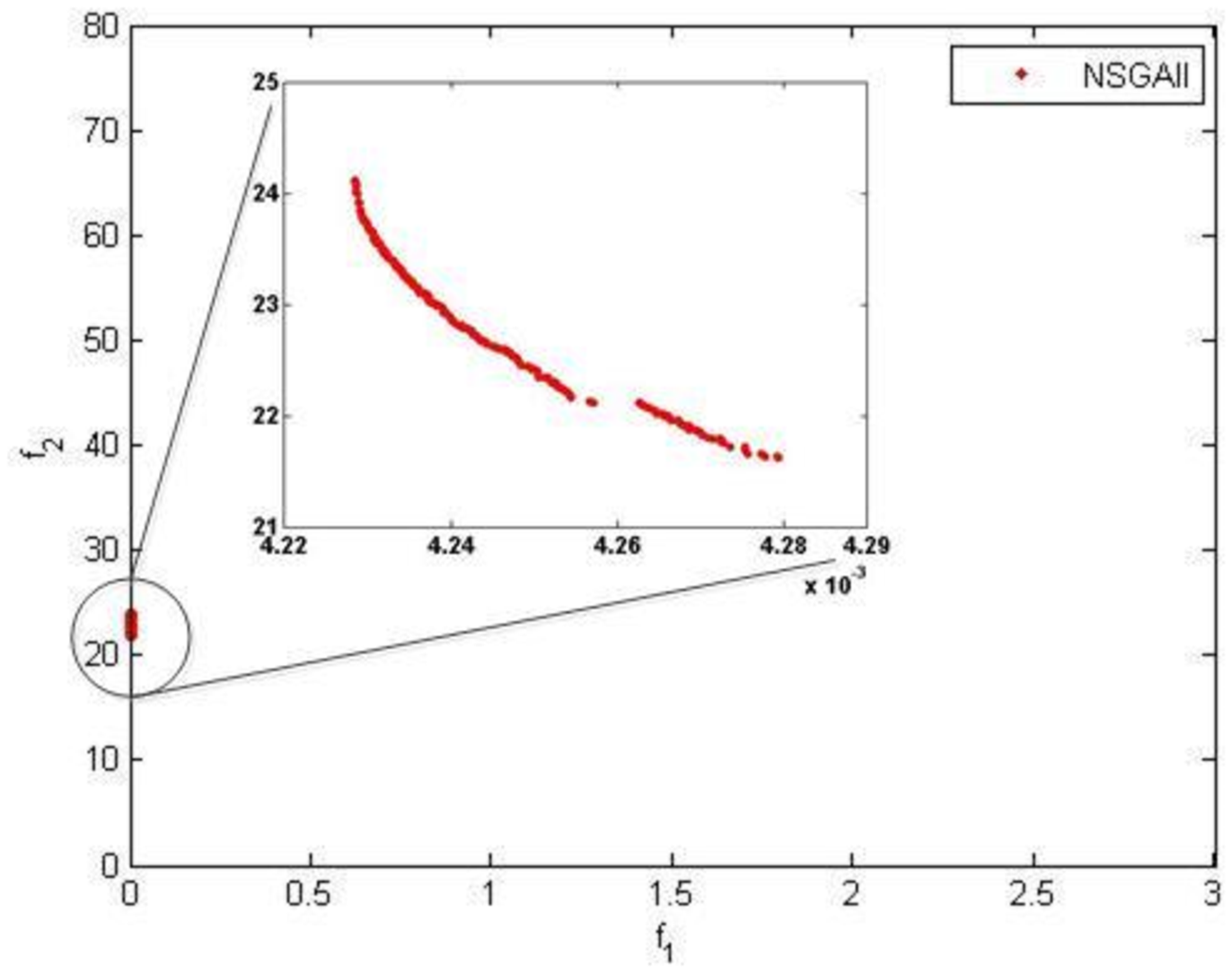,width=0.5\linewidth}\label{app:dengue:alg:a}}%
\subfigure[IBEA]{\epsfig{file=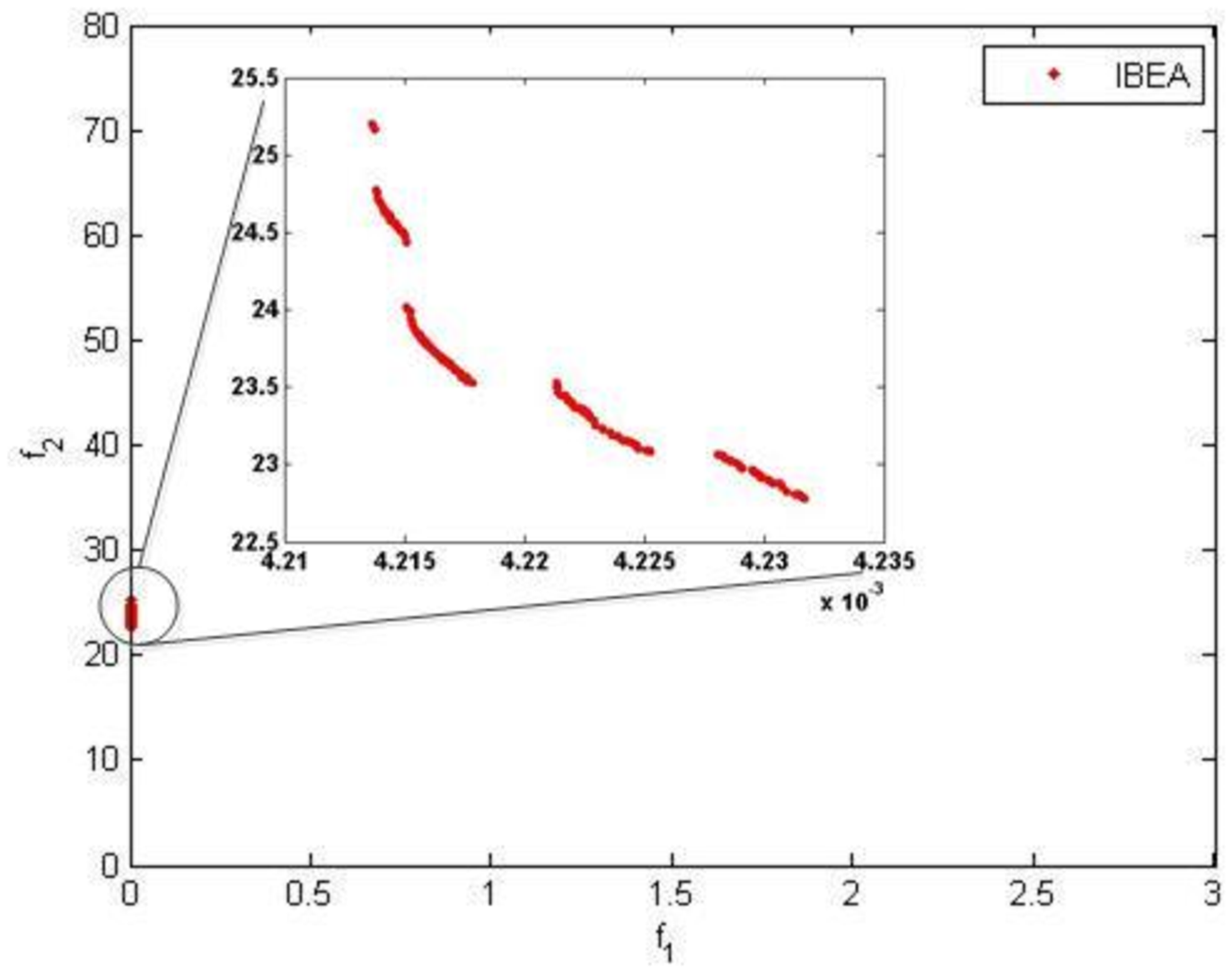,width=0.5\linewidth}\label{app:dengue:alg:b}}
\subfigure[GDE3]{\epsfig{file=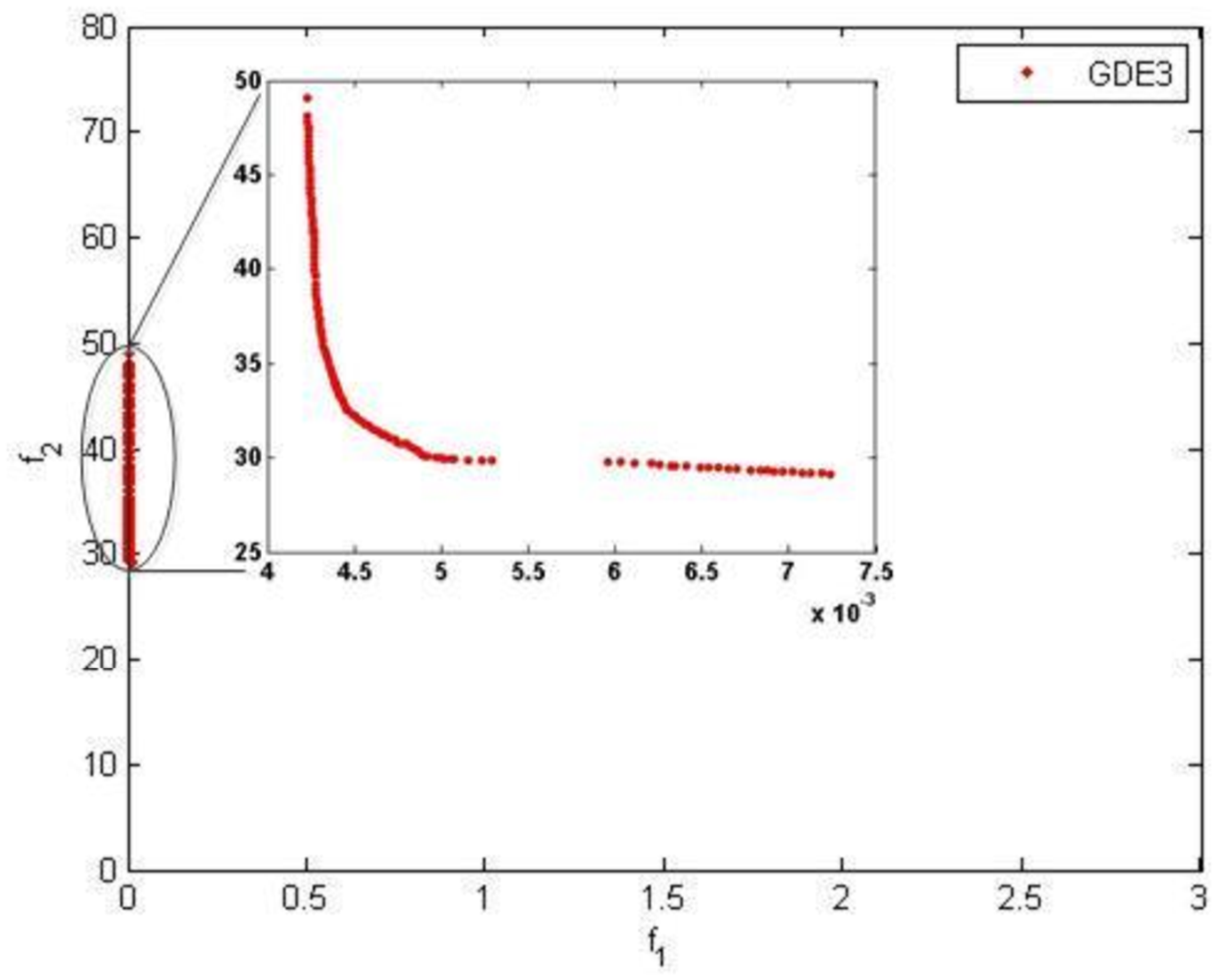,width=0.5\linewidth}\label{app:dengue:alg:c}}%
\subfigure[MOEA/D]{\epsfig{file=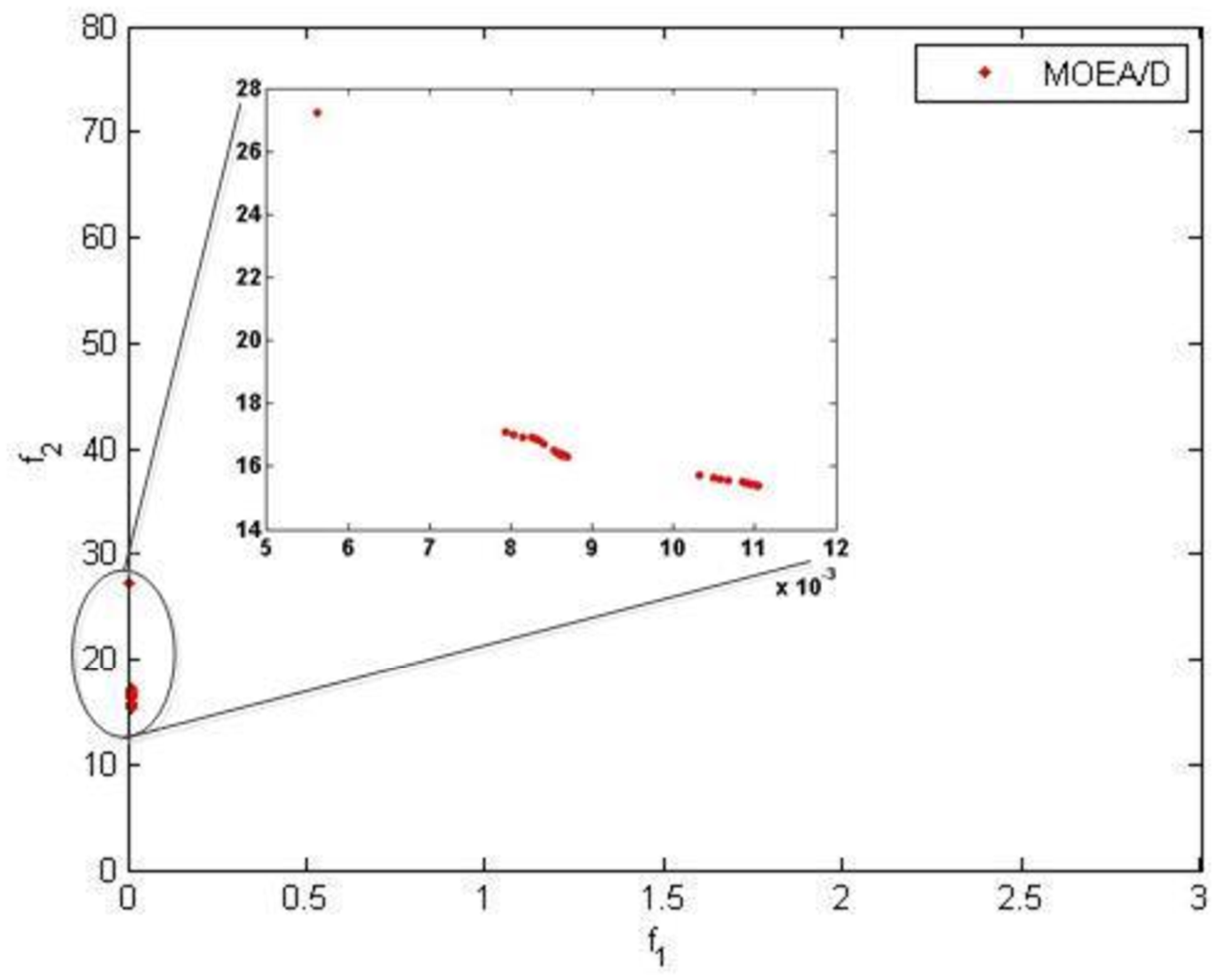,width=0.5\linewidth}\label{app:dengue:alg:d}}
\subfigure[SMPSO]{\epsfig{file=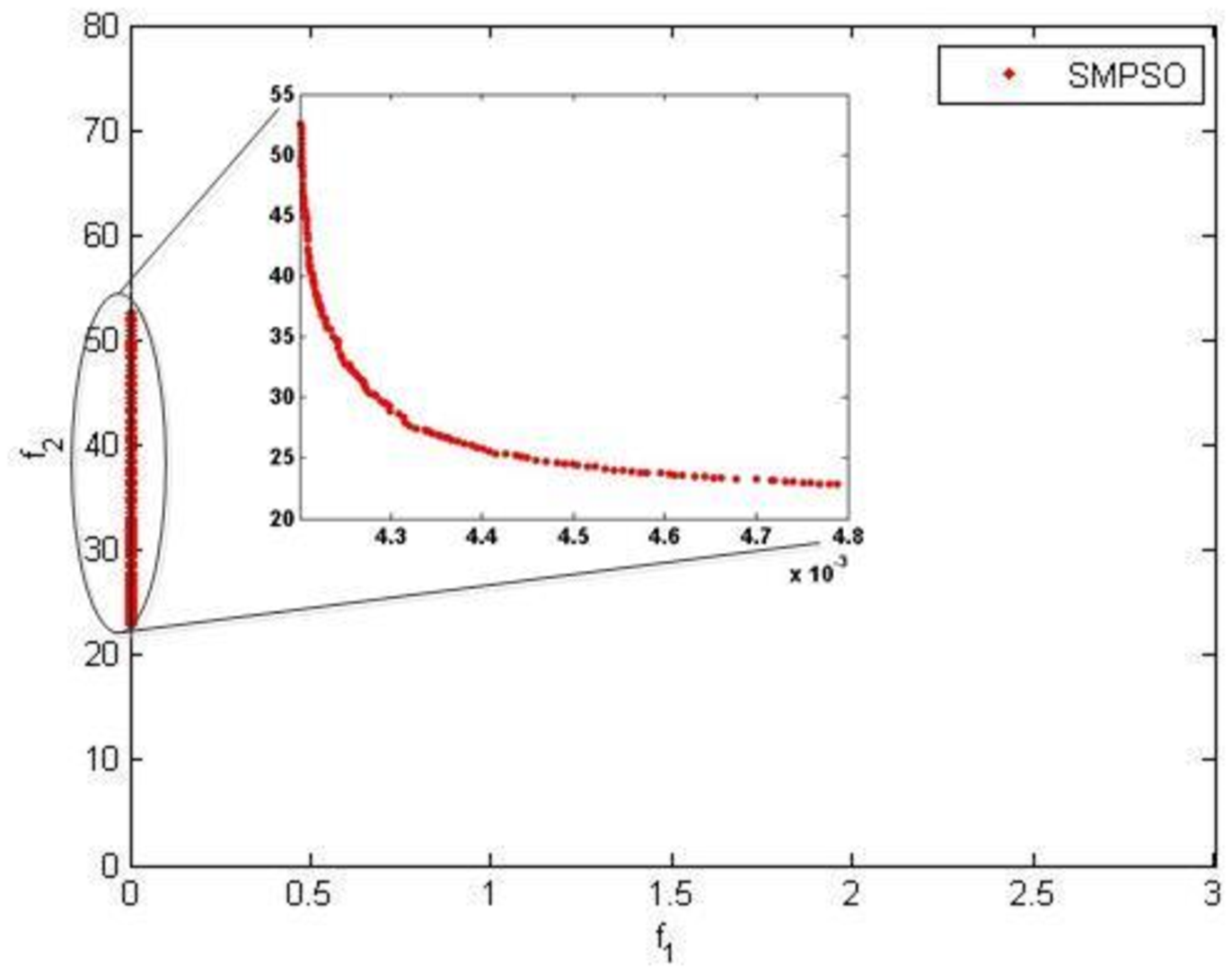,width=0.5\linewidth}\label{app:dengue:alg:e}}%
\subfigure[DDMOA2]{\epsfig{file=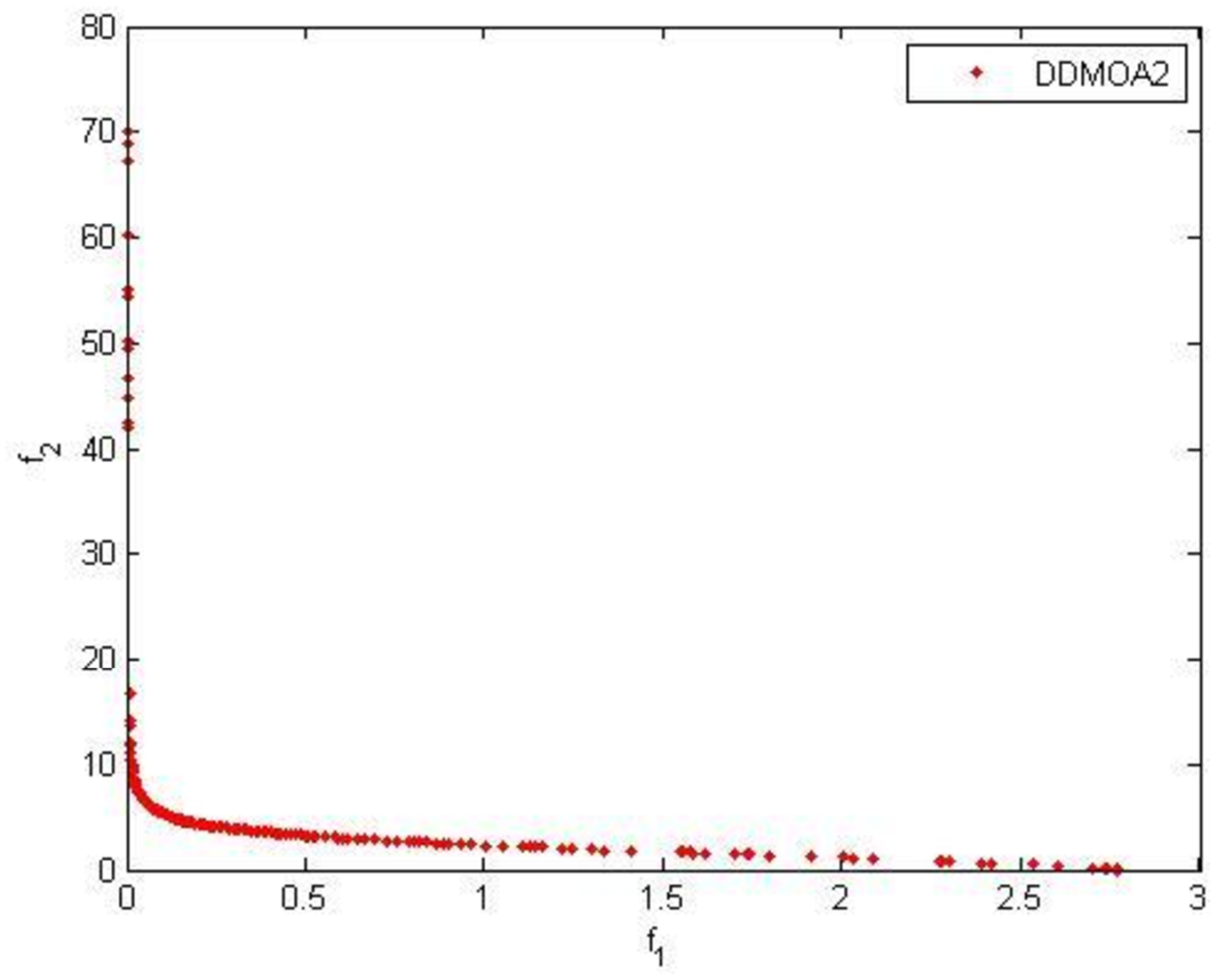,width=0.5\linewidth}\label{app:dengue:alg:f}}
\caption{Trade-off curves obtained by six different algorithms.}\label{app:dengue:alg}
\end{figure}

Figure~\ref{app:dengue:alg} depicts the sets containing all the nondominated
solutions obtained by each algorithm after 30 runs. The total infected human
population is shown in the x-axis ($f_1$). The total cost of insecticide
is shown in the y-axis ($f_2$). One can easily observe that all EMO algorithms,
with the exception of DDMOA2, face significant difficulties in obtaining a set
of well-distributed nondominated solutions in the objective space. The obtained
solutions are located in very small regions, while the majority of the search
space remains unexplored (Figures~\ref{app:dengue:alg:a}--\ref{app:dengue:alg:e}).
However, DDMOA2 is the only algorithm able to provide a good spread in the obtained
nondominated solutions (Figure~\ref{app:dengue:alg:f}). DDMOA2 extensively explores
the search space and provides a wide range of trade-off solutions. Even the visual
comparison of the obtained results allows to conclude that DDMOA2 performs
significantly better on this problem compared to all the other tested algorithms.

\begin{figure}
\centering
\subfigure[distribution of hypervolume
values]{\epsfig{file=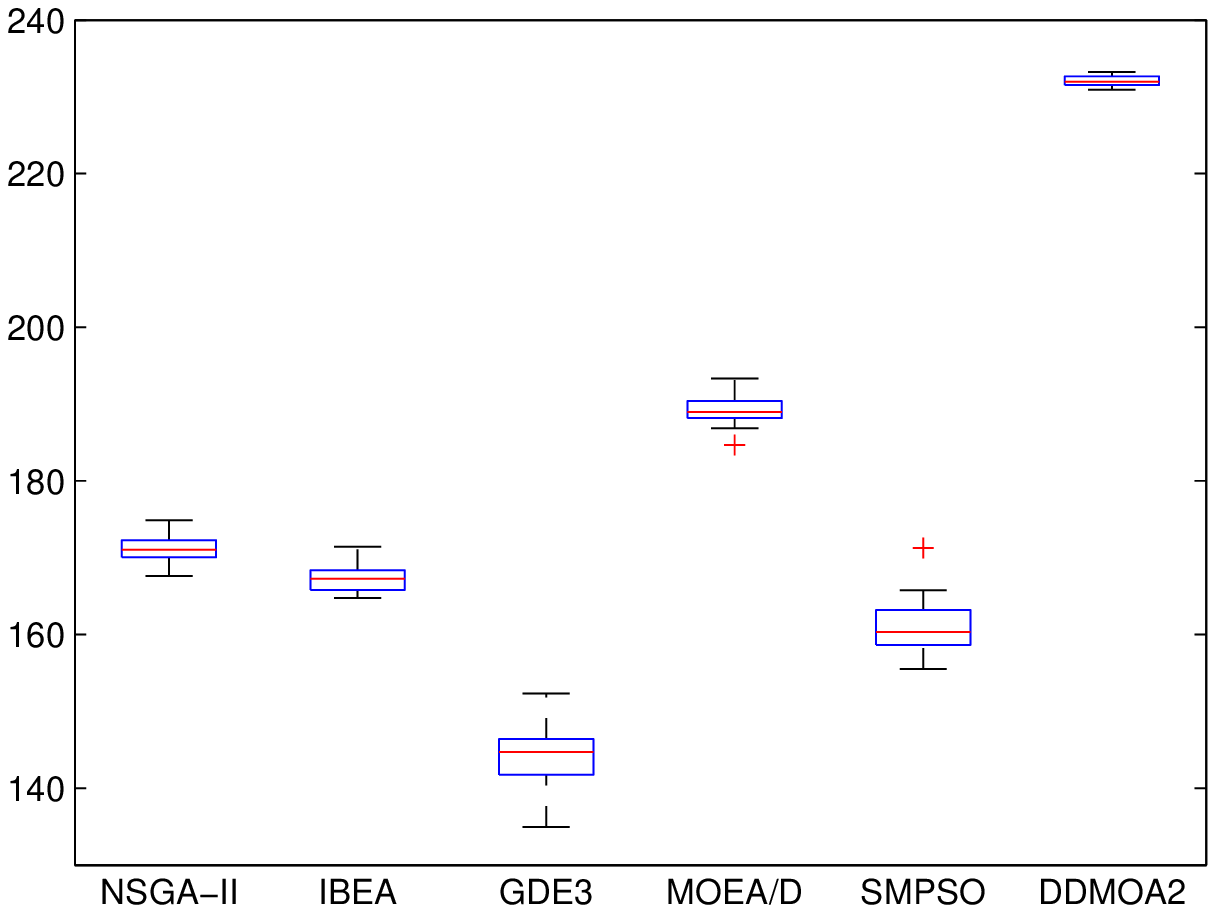,width=0.5225\linewidth}
\label{app:dengue:hv:a}}%
\subfigure[total hypervolume]{\epsfig{file=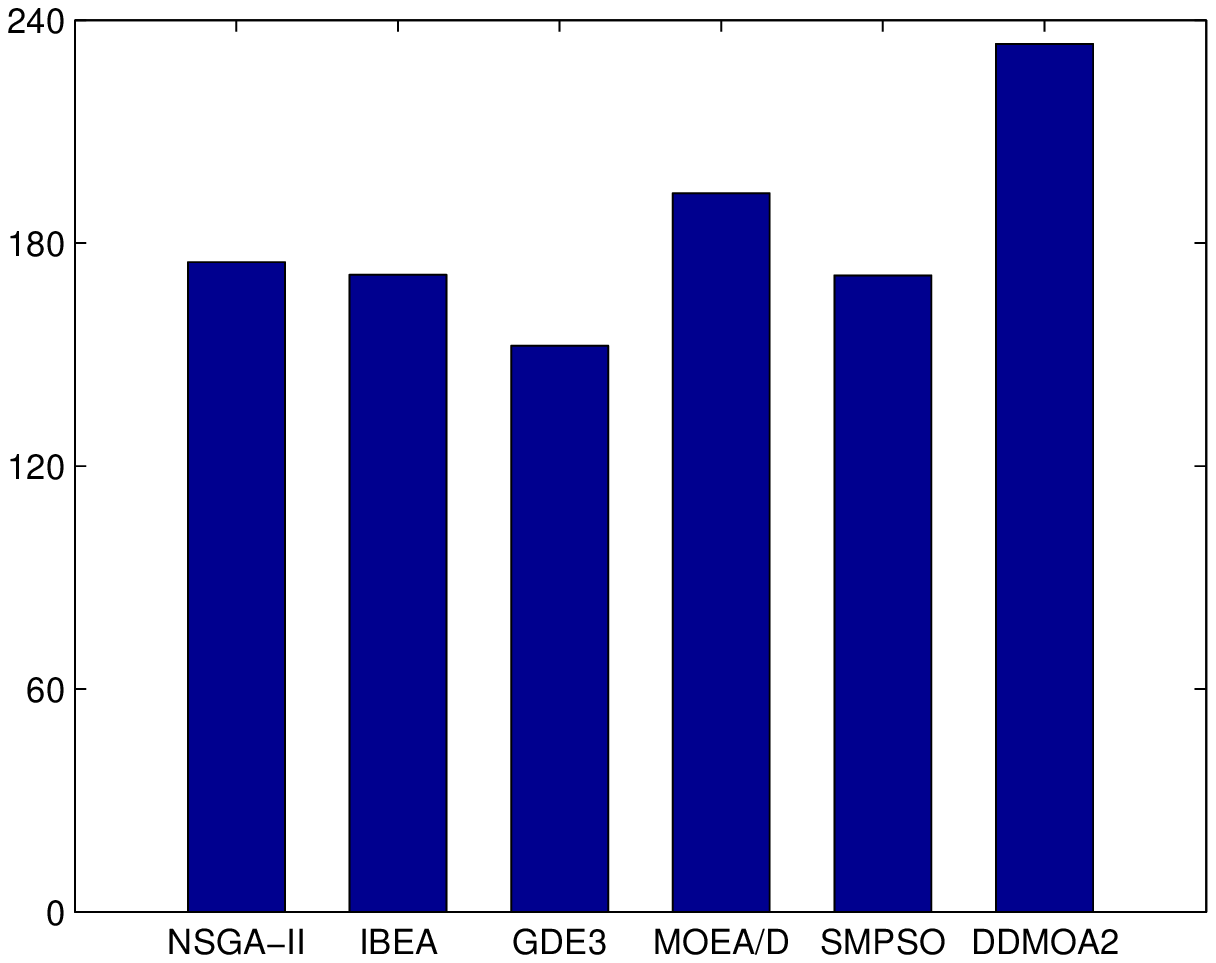,width=0.4775\linewidth}
\label{app:dengue:hv:b}}
\caption{Performance comparison of the algorithms in terms of the hypervolume
on the dengue transmission model.}\label{app:dengue:hv}
\end{figure}

To quantitatively assess the outcomes produced by the algorithms,
the hypervolume~\cite{ZitzlerThiele1998} is calculated using
$(f_1, f_2) = (3,80)$ as a reference point. Performance comparison
of the algorithms with respect to the hypervolume is shown in
Figure~\ref{app:dengue:hv}. The boxplots representing the distributions
of the hypervolume values over the runs for each algorithm are depicted
in Figure~\ref{app:dengue:hv:a}. It is interesting to note that the performance
of genetic algorithms, namely NSGA-II (dominance-based) and IBEA (indicator-based),
seems to be quite similar. However, NSGA-II gives a better spread of solutions
that results in slightly higher values of the hypervolume. Two DE-based algorithms,
namely GDE3 (dominance-based) and MOEA/D (scalarizing-based), perform differently.
GDE3 has the worst performance, while MOEA/D behaves the best without considering
DDMOA2. It can also be observed that the only PSO-based algorithm performs poorly.
In turn, DDMOA2 yields the highest values of the hypervolume when compared to
the other EMO algorithms. Moreover, the small variability of the achieved values
highlights the robustness of DDMOA2. Performance comparison with respect
to the total hypervolume achieved by all nondominated solutions obtained
by each algorithm is presented in Figure~\ref{app:dengue:hv:b}. The data
presented in this plot are consistent with the previous observations,
being DDMOA2 again the algorithm with the best performance.


\subsection{Problem Difficulties}

\begin{figure}
\centering
\includegraphics[width=0.67\textwidth]{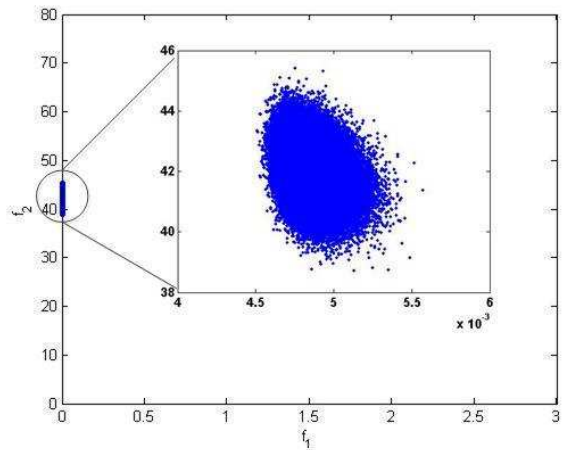}
\caption{Mapping of uniformly sampled points in the decision space into
the objective space defined by the dengue transmission model.}
\label{app:dengue:random}
\end{figure}

To better understand the difficulties in solving this multiobjective
optimization problem, $10^5$ points are sampled within the $1001$-dimensional
unit hypercube (the feasible decision space) using a uniform distribution.
Figure~\ref{app:dengue:random} illustrates the mapping of these points
into the objective space. It can be observed that the uniform distribution
of solutions in the decision space does not correspond to a uniform distribution
in the objective space. Solutions are mapped into a relatively small region
of the objective space. Thus, the probability of getting solutions in the
region shown in Figure~\ref{app:dengue:random} is much higher than anywhere
in the objective space. So there exists a pretty significant bias that
can easily deceive the operators of evolutionary algorithms.

Another difficulty that significantly exacerbates the above discussed problem
bias is related with the high dimensionality of the decision space. It is well
known that the volume of the decision space increases exponentially with the
number of decision variables. This means that in higher dimensions it is
required to maintain much larger populations in order to have a similar
coverage of the decision space as in the lower dimensions. This effect
is called \emph{the curse of dimensionality}~\cite{bellman1957}.

The discussed performance comparison of different algorithms reveals that
all of the considered state-of-the-art EMO algorithms face significant
difficulties dealing with such problem properties. Popular traditional
evolutionary operators to perform the search in the decision space (including
genetic algorithm-based, differential evolution-based and particle swarm
optimization-based variation operators) perform poorly due to aforementioned
difficulties. Moreover, most existing algorithms are typically tested on
artificial benchmarks with up to 30 decision variables. As a result, they
often become useless for large-scale problems that are highly common in
real-world applications. On the other hand, DDMOA2 is a hybrid algorithm,
which uses the concepts of traditional and stochastic optimization methods.
Its local search procedure used to find descent directions explores extensively
the decision space, working somewhat like global search in the case of this
problem. It appears to be much more effective than pure evolutionary
variation operators. Since DDMOA2 uses information about descent directions
for different objectives to produce offspring, its variation operator is an
intrinsically multiobjective one, contrary to the other algorithms. These
features become extremely useful when solving this problem, showing the
superiority of hybrid methodologies in algorithms' design.


\subsection{Dengue Dynamics}

For this particular real-world problem, the above analysis allows to conclude
that DDMOA2 works significantly better than the other considered algorithms.
However, from a decision maker's perspective, instead of comparing the
performance of EMO algorithms, the main concern is a set of optimal solutions.
Therefore, in the following we are focusing on presenting optimal solutions
corresponding to the most effective ways of fighting the dengue epidemics.

To ensure the near optimality properties of obtained trade-off solutions,
we perform local search using all the nondominated solutions returned by DDMOA2.
For this purpose, we use the normal constraint (NC) method proposed by Messac
and Mattson~\cite{Messac2004}. First, extreme solutions for all objectives
are improved using a single-objective optimization algorithm. Then,
the obtained solutions are used to construct the simplex with 300 uniformly
distributed points. Finally, for each NC subproblem, from the obtained
nondominated set we select a solution performing the best on the corresponding
NC subproblem as the initial point to optimize this subproblem. As
single-objective optimizer we use the MATLAB\textsuperscript{\textregistered}
built-in function \texttt{fmincon}. Since not all solutions optimal for the NC
subproblems may correspond to the Pareto optimal solutions, we discard all
dominated solutions from the resulting set.

\begin{figure}
\centering
\includegraphics[width=1\textwidth]{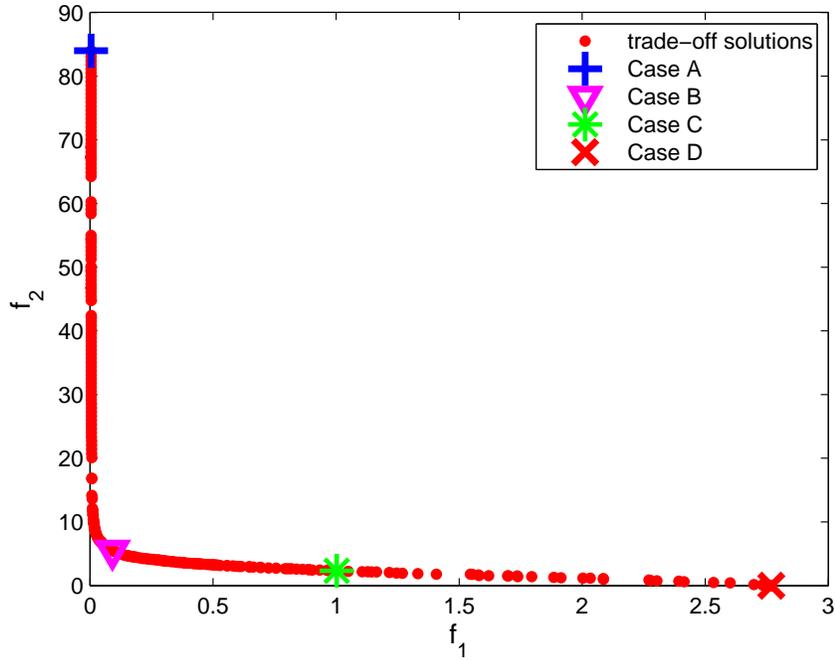}
\caption{Trade-off curve for the infected population and control.}
\label{app:dengue:all}
\end{figure}

Figure~\ref{app:dengue:all} presents all the nondominated solutions obtained
afterwards. It can be observed that the range of the trade-off curve has
been extended compared to the results shown in Figure~\ref{app:dengue:alg:f}.
Further observing this figure, one can see that for the insecticide cost
in the range $0\leq f_2\leq 5$ there is almost a linear dependency between
the infected human population ($f_1$) and the insecticide cost ($f_2$). Hence,
reducing the number of infected humans from the worst scenario to $0.5$ can be
done at a relatively low cost. However, starting from some further point,
reducing the number of infected humans can be achieved through exponential
increase in spendings for insecticide. Thus, even a small decrease in the
number of infected humans corresponds to a high increase in expenses for
insecticide. Scenarios represented by this part of the trade-off curve can be
unacceptable from the economical point of view. Furthermore, it should be noted
that even with the maximum spending it is not possible to eradicate the disease,
being 0.0042 the lowest obtained value for the infected population.

\begin{figure}
\centering
\subfigure[Case A]{\epsfig{file=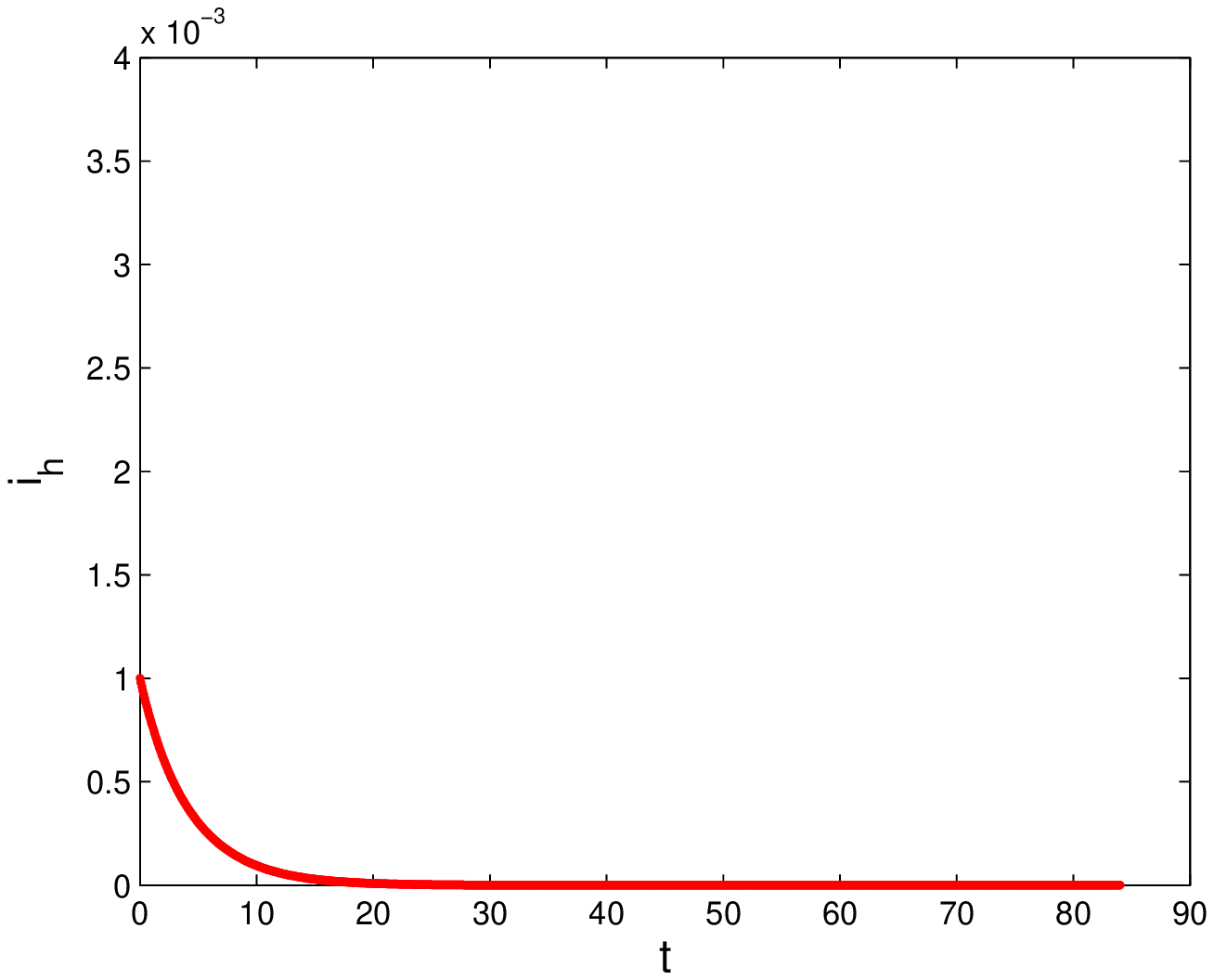,width=0.5\linewidth}\label{app:dengue:case:a}}%
\subfigure[Case B]{\epsfig{file=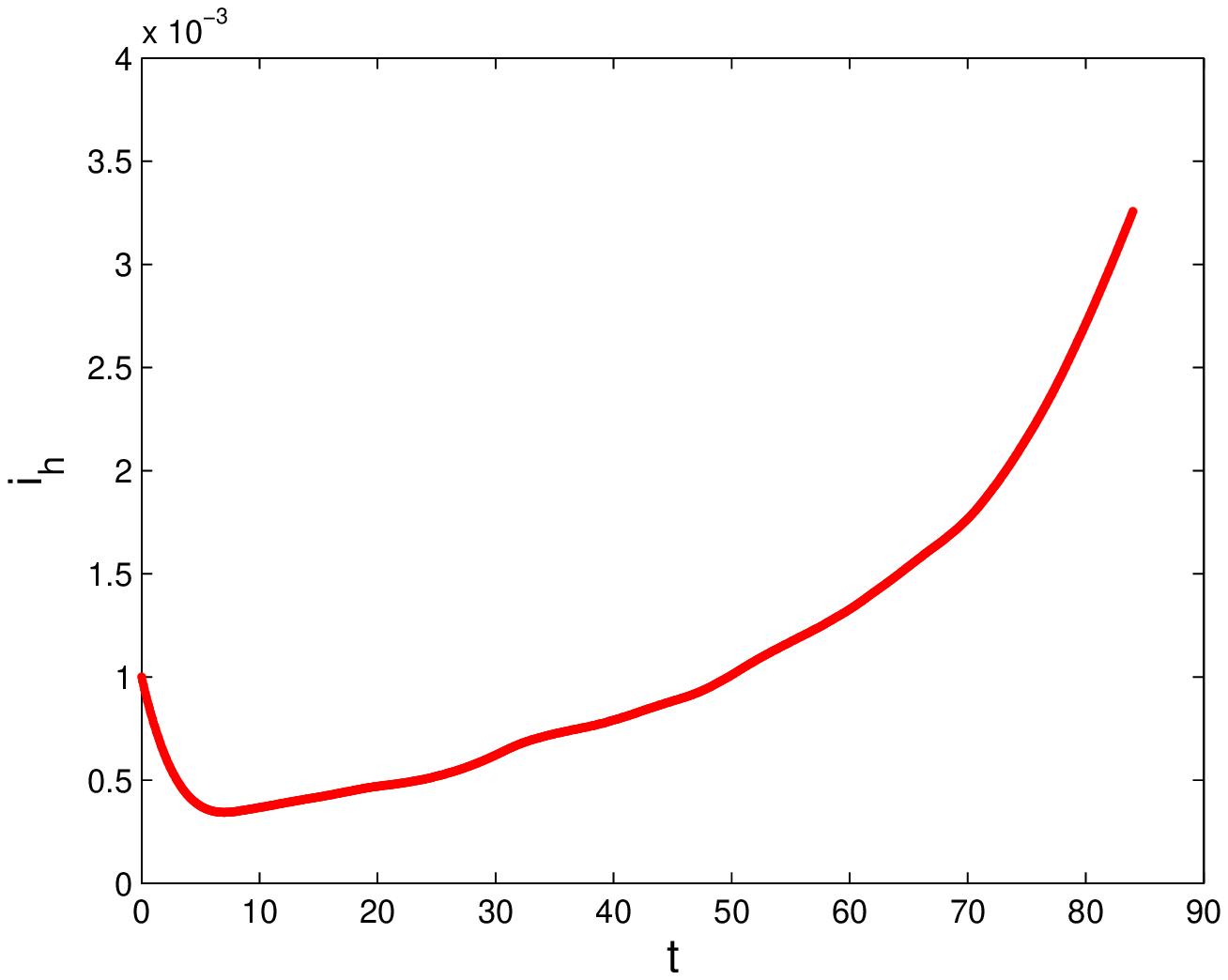,width=0.5\linewidth}\label{app:dengue:case:b}}
\subfigure[Case C]{\epsfig{file=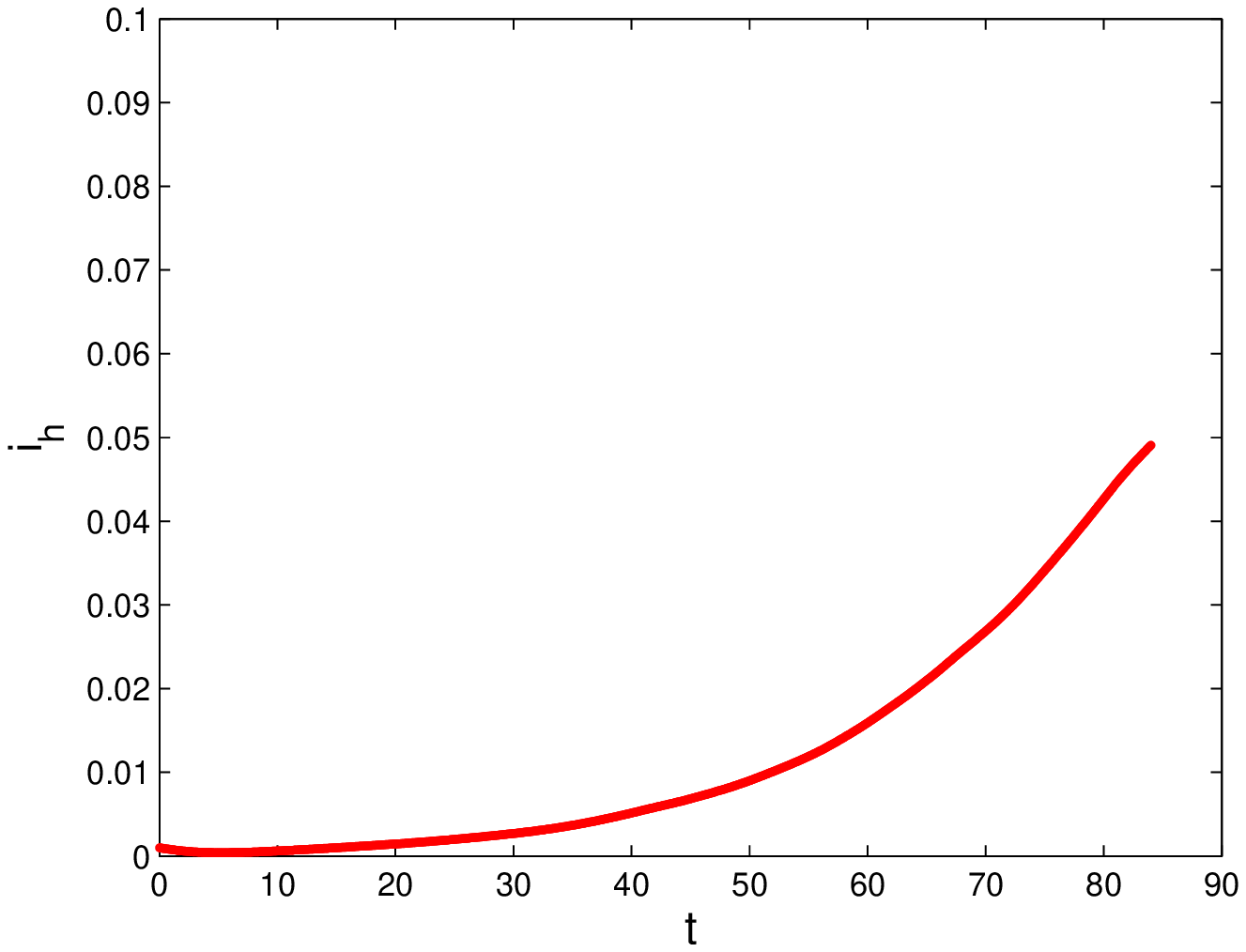,width=0.5\linewidth}\label{app:dengue:case:c}}%
\subfigure[Case D]{\epsfig{file=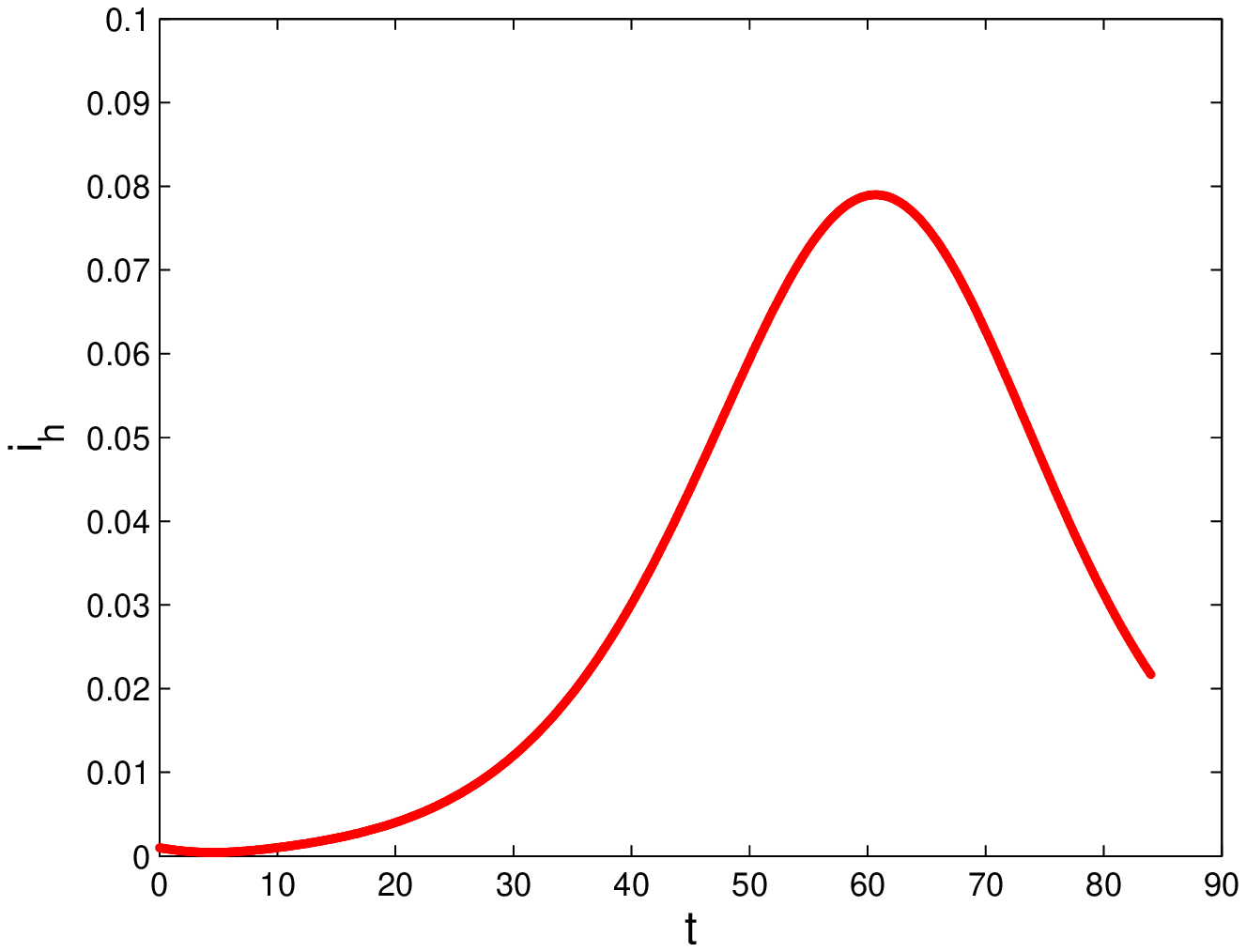,width=0.5\linewidth}\label{app:dengue:case:d}}
\caption{Different scenarios of the dengue epidemic.}\label{app:dengue:case}
\end{figure}

Additionally, Figure~\ref{app:dengue:all} presents four distinct points
(Case~A, Case~B, Case~C, and Case~D) representing different parts of the
trade-off curve, and, consequently, different scenarios of the dengue epidemic.
For each case, Figure~\ref{app:dengue:case} plots the dynamics of the numbers
of infected humans over the considered period of time.

Case A represents the medical perspective, when the number of infected humans
is the lowest. From Figure~\ref{app:dengue:case:a}, one can see that the number
of infected people decreases from the very beginning. However, it is achieved
through a huge expense for the insecticide.

Case B represents the region of the trade-off curve where an exponential
dependency between the infected humans and the control begins. From
Figure~\ref{app:dengue:case:b}, it can be seen that, at the beginning,
the number of infected humans decreases slightly. Thereafter, it grows steadily.

Case C represents the part of the trade-off curve with seemingly linear
dependency between the two objectives. Figure~\ref{app:dengue:case:c} shows
that in this case the number of infected humans grows from the very first moment.
A more rapid grow is observed in the second half of the considered period of time.

Finally, Case D represents the economical perspective, when the treatment for
infected population is neglected and the main concern is the saving from not
performing insecticide campaigns. Figure~\ref{app:dengue:case:d} illustrates
that in Case D the number of infected humans grows rapidly from the very
beginning, reaching its peak approximately on the 60th day. After that
the number of infected people decreases. This case corresponds to the
worst scenario from the medical point of view.


\subsection{Methods Comparison}

\begin{figure}
\centering
\includegraphics[width=0.8\textwidth]{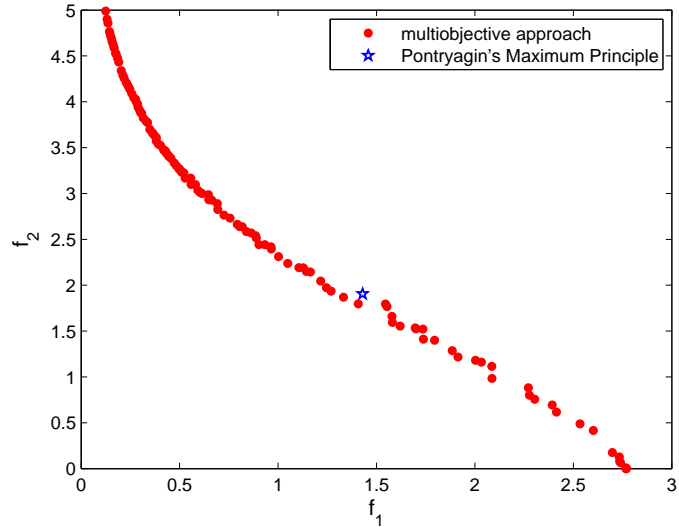}
\caption{Solutions obtained using the multiobjective approach
and Optimal Control Theory.}
\label{app:dengue:comparison}
\end{figure}

Finally, we compare the multiobjective approach discussed in this paper with
the traditional technique based on the Optimal Control Theory commonly used
to deal with such types of problems. Figure~\ref{app:dengue:comparison} shows
the region of trade-off curve with $0\leq f_2\leq 5$ and a solution obtained
using Pontryagin's Maximum Principle. From this figure, it can be seen that
there are solutions on the trade-off curve obtained using the multiobjective
approach, which dominate the solution obtained using the Optimal Control Theory.
Additionally, as it is shown in Figure~\ref{app:dengue:all} and can be seen
in Figure~\ref{app:dengue:comparison}, the multiobjective approach allows
to find a number of other solutions that are incomparable with respect to the
solution obtained using Pontryagin's Maximum Principle. All these trade-off
solutions provide more valuable information about the dengue epidemic,
offering a variety of potential strategies to fight the disease. This highlights
the advantages of a multiobjective approach compared to traditional techniques
for finding the optimal control.


\section{Conclusions}
\label{sec:5}

We proposed a multiobjective approach to find the optimal control
to manage financial expenses caused by the outbreak of the dengue epidemic.
We formulated the optimization problem with two objectives, including the system
of differential equations modelling the disease. The first objective represents
expenses due to the infected population. The second objective represents
the cost of applying insecticide in order to fight the disease. We sought
the optimal control simultaneously optimizing these clearly conflicting objectives.
The obtained trade-off solutions reveal different perspectives on applying insecticide:
a low number of infected humans can be achieved spending larger financial recourses,
whereas low spendings for prevention campaigns result in significant portions
of the population affected by the disease. At the same time, a number of other
solutions represent different trade-offs between the objectives. Once the whole
range of optimal solutions is obtained, the final decision on the control strategy
can be made taking into consideration the available financial resources and goals
of public health care. The analysis of different approaches to find the optimal
control in the proposed model shows that the multiobjective approach presents
clear advantages compared to the traditional approach based
on the Optimal Control Theory. Additionally, the performance comparison
of different multiobjective optimization algorithms is carried out on the
optimization problem resulting from the proposed modelling approach. The obtained
results show that the hybrid methodology significantly outperforms the traditional
evolutionary algorithms. DDMOA2 appears to be the only among considered algorithms
capable to find a wide range of trade-off solutions, whereas the other tested
EMO algorithms face significant difficulties in solving this problem:
the obtained solutions are located in small regions of the objective space.
We believe that the promising results presented here will promote
multiobjective modelling and will motivate further research
in the design of hybrid and local search-based multiobjective optimization
algorithms, which until now received relatively limited attention.


\subsection*{Acknowledgments}

This work has been supported by the \emph{Portuguese Foundation for Science
and Technology} (FCT) in the scope of projects UID/CEC/00319/2013
(ALGORITMI R\&D Center) and UID/MAT/04106/2013 (CIDMA).



\end{document}